\magnification=\magstep1

\def\refs{\medskip\hangindent=25pt\hangafter=1\noindent}

\centerline {ON POSSIBLE MIXING RATES FOR SOME STRONG
MIXING CONDITIONS} 
\centerline {FOR N-TUPLEWISE INDEPENDENT
RANDOM FIELDS}
\bigskip

\noindent Richard C.\ Bradley \hfil\break
Department of Mathematics \hfil\break
Indiana University \hfil\break
Bloomington \hfil\break 
Indiana 47405 \hfil\break

\noindent E-mail address: bradleyr@indiana.edu 
\hfil\break
\vskip 0.5 in

   {\bf Abstract.} For a given pair of positive
integers $d$ and $N$ with $N \geq 2$,  
for strictly stationary random
fields that are indexed by the $d$-dimensional 
integer lattice and satisfy $N$-tuplewise
independence, the dependence coefficients
associated with the $\rho$-, $\rho'$-, and $\rho^*$-mixing 
conditions can decay together at an arbitrary rate.
Another, closely related result is also established. 
In particular, these constructions provide classes of
examples pertinent to limit theory for random fields that
involve such mixing conditions together with certain
types of ``extra'' assumptions on the
marginal and bivariate (or $N$-variate) distributions.
\hfil\break     

\hfil\break

\hfil\break

\vfill\eject
\noindent {\bf 1. Introduction}
\hfil\break

   In certain types of limit theorems for random fields
under strong mixing conditions, there are sometimes 
``extra conditions'' involving the marginal and (say)
bivariate distributions.
For example, in a paper by the author and Tran [7], 
a central limit theorem for some kernel-type estimators 
of probability density was proved for strictly
stationary random fields that
satisfy certain strong mixing conditions as well as
certain ``extra'' assumptions:
an absolutely continuous marginal distribution
with a continuous density function, and a
dependence condition on the bivariate distributions. 
In a subsequent paper, the author [3] constructed a
class of random fields that satisfy all of those
conditions (including the ``extra'' ones);
and by altering the choices of parameters in those 
examples, one can get practically arbitrary ``mixing 
rates'' within certain rather loose constraints.
Later, Tone [24, Theorem 5.2] proved a functional 
central limit
theorem for empirical processes from strictly stationary
$\rho'$-mixing random fields that satisfy the following
``extra'' conditions: the marginal distribution is
absolutely continuous and (for convenience) supported
on the unit interval, and there exists a positive 
constant $C$ such that the bivariate distributions are absolutely continuous with bivariate density functions
that are (at any given point $(x,y)$) bounded above by
$C$ times the product of the marginal densities. 
(The mixing conditions in the papers [7] and [3]
cited above were quite different from $\rho'$-mixing,
and will not be dealt with further in this note.
The $\rho'$-mixing condition itself will be defined 
below.)\ \ 
The author [5] constructed a  
class of strictly stationary $\rho'$-mixing random fields
with arbitrary mixing rates.
That construction gave further information on other
related mixing conditions (more on
that below);
however, it did not address the question whether such 
$\rho'$-mixing examples exist which also satisfy 
certain types  of ``extra conditions'' (on the marginal and bivariate distributions) such as in the result of
Tone [24] cited above. 
In this note, that question will be answered
affirmatively, \`a la the paper [3] cited above.
\hfil\break 

   Those ``extra conditions'' in the cited result 
from [24] are
trivially satisfied in the case where the random 
variables are (say) uniformly distributed on the unit 
interval and are pairwise independent.
In the examples that are constructed in this note,
the $\rho'$-mixing condition will 
(for appropriate choices of parameters)
be satisfied with an
arbitrary mixing rate, the marginal distributions
will be uniform on the unit interval, and the random
variables will be pairwise independent (and
satisfy an even stronger ``independence'' property).
The examples given below will also involve
three other mixing conditions closely related to
$\rho'$-mixing; motivations for that will be given
below.  \hfil\break

   Before the main results are stated, some notations
and definitions need to be given.     

Throughout this paper, the setting will be a
probability space $(\Omega, {\cal F}, P)$ rich enough
to accommodate all random variables defined here.
For any collection $(Y_i, i \in I)$  of random
variables defined on this probability space, the
$\sigma$-field of events generated by this collection
will be denoted $\sigma (Y_i, i \in I)$.

   Let ${\bf N}$ (resp.\ ${\bf R}$ resp.\ ${\bf Z}$)
denote the set of all positive integers (resp.\ 
all real numbers resp.\ all integers).

   Suppose $d$ is a positive integer.
A given element $k \in {\bf Z}^d$ will be represented by
$k := (k_1, k_2, \dots, k_d)$.
For any element $k \in {\bf Z}^d$, define the
Euclidean norm 
$\|k\| := (k_1^2 + k_2^2 + \dots + k_d^2)^{1/2}$.
The origin in ${\bf Z}^d$ will be denoted with boldface 
${\bf 0} := (0,0,\dots, 0)$.
For any two nonempty disjoint subsets $S$ and $T$
of ${\bf Z}^d$, define the positive number
${\rm dist}(S,T) := \min_{s \in S, t \in T} \|s-t\|$.

   (Of course if $d = 1$, then
$\|k\| = |k|$, ${\bf 0} = 0$, and 
${\rm dist}(S,T) = \min_{s \in S, t \in T} |s-t|$.)
\hfil\break  

{\bf Definition 1.1.}  
Suppose $d$ is a positive integer,
and $X := (X_k, k \in {\bf Z}^d)$ is a (not necessarily
stationary) random field on $(\Omega, {\cal F}, P)$.
Suppose $N \geq 2$ is an integer.  
The random field $X$ satisfies 
``$N$-tuplewise independence'' if the following holds:  
For every choice of $N$ distinct elements 
$k(1),\ k(2),\ \dots,\ k(N) \in {\bf Z}^d$, the
random variables
$X_{k(1)},\ X_{k(2)},\ \dots,\ X_{k(N)}$ are independent.
\hfil\break

   The notion of $N$-tuplewise independence --- and
in particular, pairwise independence (the case $N = 2$)
--- is of interest in its own right.
Limit theory for sequences of pairwise independent 
random variables has been developed in a consider
number of papers, including the ones by Etemadi [12][13]
and Etemadi and Lenzhen [14].
The paper [13] cites several references for practical
applications of pairwise independence 
(e.g.\ in statistics and computer science).
Also, some limitations of that theory have been 
illustrated with, for example, nontrivial pairwise 
independent counterexamples to the central limit 
theorem, such as constructions in the papers 
of Janson [17], of Cuesta and Matr\'an [9], and
of the author and Pruss [6].
The example in the latter paper satisfies 
$N$-tuplewise independence, where $N$ is an
arbitrary integer $\geq 2$ specified beforehand.
In the book by McWilliams and Sloane [19], in
connection with the construction of 
error-correcting codes, some methodology
involves the construction of ``big'' $N$-tuplewise 
independent random vectors from ``small'' ones.
The paper by Rosenblatt and Slepian [23] dealt with the question of whether (for a given $N$ and $n$) 
$N$-tuplewise independence can hold for $n$th order Markov chains with certain restrictions on the state space.
\hfil\break

   Now let us turn our attention to mixing conditions.
\hfil\break   

{\bf Definition 1.2.}
For any two $\sigma$-fields 
${\cal A}$ and ${\cal B} \subset {\cal F}$, define the
following two measures of dependence:
$$ \eqalignno {\alpha ({\cal A}, {\cal B}) &:=
\sup_{A \in {\cal A}, B \in {\cal B}}
|P(A \cap B) - P(A)P(B)| \quad {\rm and} & (1.1) \cr
\rho({\cal A}, {\cal B}) &:= 
\sup |{\rm Corr}(f,g)|     & (1.2) \cr}$$
where the latter supremum is taken over all pairs of
square-integrable random variables $f$ and $g$ such that
$f$ is ${\cal A}$-measurable and $g$ is 
${\cal B}$-measurable.
The quantity $\rho({\cal A}, {\cal B})$
is known as the ``maximal correlation coefficient''
between ${\cal A}$ and ${\cal B}$.
The following inequality is elementary
(see e.g.\ [4, v1, Proposition 3.11(b)]):  
For any two $\sigma$-fields ${\cal A}$ and ${\cal B} \subset {\cal F}$,
$$ 4 \alpha ({\cal A}, {\cal B}) \leq
\rho ({\cal A}, {\cal B}).  \eqno (1.3) $$   
       
   {\bf Definition 1.3.} Again suppose $d$ is a positive
integer, and
$X := (X_k, k \in {\bf Z}^d)$ is a (not 
necessarily stationary) random field on 
$(\Omega, {\cal F}, P)$.
(No assumption of ``$N$-tuplewise independence.'')\ \
For each positive integer $n$, 
referring to (1.1) and (1.2), define the following
four ``dependence coefficients'':  First,
$$ \eqalignno{
\alpha(n) = \alpha(X,n) &:=
\sup \alpha (\sigma(X_k, k \in S), \sigma(X_k, k \in T))
\quad {\rm and}
   & (1.4) \cr  
\rho(n) = \rho(X,n) &:=
\sup \rho (\sigma(X_k, k \in S), \sigma(X_k, k \in T)) 
   & (1.5) \cr}$$
where in each of (1.4) and (1.5) the supremum is taken 
over all pairs of
sets $S$ and $T \subset {\bf Z}^d$ such that
for some integer $j$ and some $u \in \{1, 2, \dots, d\}$,
one has that
$S = \{k \in {\bf Z}^d: k_u \leq j\}$ and
$T = \{k \in {\bf Z}^d: k_u \geq j+n\}$.
Next, 
$$ \rho'(n) = \rho'(X,n) :=
\sup \rho (\sigma(X_k, k \in S), \sigma(X_k, k \in T)) 
\eqno (1.6) $$
where this supremum is taken over all pairs of
sets $S$ and $T \subset {\bf Z}^d$ such that for
some $u \in \{1, 2, \dots, d\}$ and some pair of
nonempty, disjoint sets $G$ and $H \subset {\bf Z}$
such that 
${\rm dist}(G,H) \geq n$,
one has that   
$S = \{k \in {\bf Z}^d: k_u \in G\}$ and
$T = \{k \in {\bf Z}^d: k_u \in H\}$.
(Note that the sets $G$ and $H$ may
be ``interlaced,'' with each one containing elements
between ones in the other set.)\ \ Finally,
$$ \rho^*(n) = \rho^*(X,n) :=
\sup \rho (\sigma(X_k, k \in S), \sigma(X_k, k \in T)) 
\eqno (1.7) $$
where this supremum is taken over all pairs of
nonempty, disjoint
sets $S$ and $T \subset {\bf Z}^d$ such that
${\rm dist}(S,T) \geq n$.     

   Obviously, for each $n \in {\bf N}$, by (1.3),
$$ 0 \leq 4\alpha(n) \leq \rho(n) \leq \rho'(n) 
\leq \rho^*(n) \leq 1.
\eqno (1.8)$$
Also obviously, the sequences 
$(\alpha(n), n \in {\bf N})$,
$(\rho(n), n \in {\bf N})$,
$(\rho'(n), n \in {\bf N})$, and
$(\rho^*(n), n \in {\bf N})$ are each nonincreasing as
$n$ increases.
The random field $X$ is said to be \hfil\break
``strongly mixing,'' or ``$\alpha$-mixing,'' 
if $\alpha(n) \to 0$ as $n \to \infty$, \hfil\break
``$\rho$-mixing'' if $\rho(n) \to 0$ as $n \to \infty$, 
\hfil\break
``$\rho'$-mixing'' if $\rho'(n) \to 0$ as $n \to \infty$, and
\hfil\break
``$\rho^*$-mixing'' if $\rho^*(n) \to 0$ as $n \to \infty$.
\hfil\break

   Our focus will be on strictly stationary
random fields.
There is a vast literature on limit theory for
strictly stationary random fields satisfying
the mixing conditions above and closely related ones. 
See e.g.\ [4, v3], [7], [15], [16], [21], 
[22], [24], and [25], to name just a few references. 
As in the well known paper of Kesten and O'Brien [18]
(involving random sequences --- the case $d=1$), there has
been an ongoing interest in the question of what 
``mixing rates'' are possible for the 
above mixing conditions and related ones, especially
under strict stationarity.
The papers [1], [3], and [5] deal with that question
for various mixing conditions for strictly
stationary random fields.
Theorems 1.4 and 1.5 below will give some results
on the question of what ``mixing rates'' can occur
for such random fields that are also $N$-tuplewise 
independent.
\hfil\break

   A couple of background facts are in order.
First (even without stationarity), for the case $d=1$ 
(random sequences) one trivially has that 
$\rho^*(n) = \rho'(n)$.  
Second, if the random field $X$ is strictly stationary, then 
(i) for the case $d \geq 2$, one has that
$4 \alpha(n) \leq \rho(n) \leq 2\pi \alpha(n)$, and 
(ii) for any positive integer $d$ (including $d = 1$),
analogous inequalities hold 
(in comparison to $\rho'(n)$ and $\rho^*(n)$ respectively)
for dependence coefficients
$\alpha'(n)$ and $\alpha^*(n)$ that are defined
\`a la (1.6) and (1.7) using the measure of dependence
$\alpha ({\cal A}, {\cal B})$.
(See [4, v3, Theorem 29.12].)\ \ 
Thus there is no real reason to formally 
define such dependence coefficients 
$\alpha'(n)$ and $\alpha^*(n)$; and the treatment of
the dependence coefficient $\alpha(n)$ (separately
from $\rho(n)$) is meaningful only in the case
$d = 1$. 
\hfil\break

   We shall now state the two main results of this note, 
and then give some further motivations for them.
\hfil\break

   {\bf Theorem 1.4.}  {\sl Suppose $d$ and $N$ are positive 
integers such that $N \geq 2$.
Suppose $(c_1, c_2, c_3, \dots)$ is a nonincreasing 
sequence of numbers in the closed unit interval $[0,1]$.
Then there exists a strictly stationary random field
$X := (X_k, k \in {\bf Z}^d)$ with the following
properties:

   (A) The random variable $X_{\bf 0}$ is uniformly
distributed on the interval $[0,1]$.

   (B) The random field $X$ satisfies $N$-tuplewise
independence.

   (C) For each $n \in {\bf N}$, 
$4\alpha(n) = \rho(n) = \rho'(n) = \rho^*(n) = c_n$.}
\hfil\break

{\bf Theorem 1.5.}  {\sl Suppose $d$ and $N$ are 
integers such that $d \geq 2$ and $N \geq 2$.
Suppose $(c_1, c_2, c_3, \dots)$ is a nonincreasing 
sequence of numbers in the closed unit interval $[0,1]$.
Then there exists a strictly stationary random field
$X := (X_k, k \in {\bf Z}^d)$ with the following
properties:

   (A) The random variable $X_{\bf 0}$ is uniformly
distributed on the interval $[0,1]$.

   (B) The random field $X$ satisfies $N$-tuplewise
independence.

   (C) For each $n \in {\bf N}$, $\rho^*(n) = 1$;
also $\rho(1) = \rho'(1) = 1$. 
   
   (D) For each $n \geq 2$, 
$\rho(n) = \rho'(n) = c_n$.}
\hfil\break

   In these two theorems, there is no assumption that 
$c_n \to 0$ as $n \to \infty$.
For some limit theory for random fields under just
the ``partial mixing'' condition
$\lim_{n \to \infty} \rho^*(n) < 1$, see e.g.\  
Peligrad and Gut [21].
\hfil\break 

   The constructions in the proofs of Theorems 1.4 and 1.5 
are based on a
well known simple class of random vectors (given in 
Definition 2.5 in Section 2) with very nice 
properties related to both $N$-tuplewise independence 
and maximal correlation coefficients.
However, this class does not seem to provide the
leverage for (i) ``separating'' the dependence
coefficients $\rho(n)$ and $\rho'(n)$ (as was
done in the construction in [5] alluded to above) or 
(ii) getting a value of $\rho(1)$ (or $\rho'(1)$)
that is less than that of $\rho^*(1)$.
In Theorem 1.5, the explicit mention of the equations 
$\rho(1) = \rho'(1) = 1$ is intended only as a reminder 
of that latter fact.
\hfil\break   
  
   In [4, v3, Theorem 26.8 and its subsequent Note 3] 
it was shown that for
strictly stationary random sequences (the case $d=1$)
with a purely non-atomic marginal distribution,
the simultaneously behavior of the dependence 
coefficients $\alpha(n)$, $\rho(n)$, and 
$\rho^*(n)$ could be practically arbitrary subject to
(1.8) and the sentence after it.
(In that ``subsequent Note 3,'' the word ``atomic''
should be ``nonatomic.'')\ \ 
In [5, Theorem 1.9], it was shown that for strictly
stationary random fields in the case $d \geq 2$ with
a purely non-atomic marginal distribution, 
an analogous fact holds for $\rho(n)$, $\rho'(n)$,
and $\rho^*(n)$.
However, those two constructions did not address
the question of whether such examples exist in which 
the bivariate distributions could have nice ``extra''
properties such as the ones assumed in the result of 
Tone [24] cited above.
To a certain extent, Theorems 1.4 and 1.5
address that question.  
Theorem 1.4 shows that for general $d \geq 1$,
the cited result in [24],
which involves $\rho'$-mixing with a
mixing rate that can be arbitrarily slow along with
the ``extra'' assumptions alluded to above on the 
marginal and 
bivariate distributions, cannot be obtained from any
similar result involving $\rho$-mixing with a 
``mixing-rate'' assumption (such as the ``logarithmic''
mixing-rate assumption 
$\sum_{n=1}^\infty \rho(2^n) < \infty$
that plays a key role in much of the limit theory under 
$\rho$-mixing). 
For $d=1$ it also shows that any similar result 
(that includes those ``extra'' assumptions) under 
$\rho$-mixing (or even $\rho^*$-mixing)
with a particular mixing-rate assumption, 
cannot be obtained from an analogous result involving 
$\alpha$-mixing with a faster mixing rate.
Theorem 1.5 shows that for the case $d \geq 2$, 
the cited result in [24] also cannot be obtained 
from any similar result under $\rho^*$-mixing or even 
under assumptions that include $\rho^*(n) < 1$ for some 
$n \geq 1$.
\hfil\break

   Theorem 1.4 will be proved in Section 3 below, 
after some preliminary work in Section 2.
Then Theorem 1.5 will be proved in Section 4
(with Theorem 1.4 itself serving as a key 
``building block'' in that proof).
\hfil\break

\noindent {\bf 2. Preliminaries} \hfil\break

   First here are some notations that will be used
in the rest of this paper.

   For a given pair of nonempty sets $A$ and $B$, the 
notation $B^A$ will of course mean the set of all mappings 
from $A$ to $B$.
A notation such as $(r_s, s \in S)$ will refer to a
particular mapping from a given nonempty set $S$ to ${\bf R}$
(the one in which $s \mapsto r_s$ for each $s \in S$).

   When a notation of the form $a^b$ appears in a
subscript, superscript, or exponent, it will be
written $a \uparrow b$.

   For any given set $S \in {\bf Z}^d$ and any
$v \in {\bf Z}^d$, the notation $S+v$ will of course
refer to the set of elements of the form $s+v$
for $s \in S$.

   The cardinality of a given set $S$ will be denoted
${\rm card}\thinspace S$. 
\hfil\break

   Now in order to avoid some clutter in
arguments later on,
we shall give a review of three standard techniques for
constructing strictly stationary random fields from other
(e.g.\ ``building block'') random fields.
\hfil\break

   {\bf Remark 2.1.} (A) Suppose $d$ is a positive integer.
Suppose that for each $n \in {\bf N}$, 
$Y^{(n)} := (Y^{(n)}_k, k \in {\bf Z}^d)$ is a 
strictly stationary random field.
Suppose these random fields $Y^{(n)}, n \in {\bf N}$
are independent of each other.
Suppose $f:{\bf R}^{\bf N} \to {\bf R}$ is a Borel
function, and the random field 
$X := (X_k, k \in {\bf Z}^d)$ is defined by
$X_k := f(Y^{(1)}_k, Y^{(2)}_k, Y^{(3)}_k, \dots)$
for $k \in {\bf Z}^d$.
Then the random field $X$ is strictly stationary.

   (B) In the context of (A) above, if also $N \geq 2$
is an integer and each of the random fields
$Y^{(n)}, n \in {\bf N}$ satisfies $N$-tuplewise
independence, then the random field $X$ satisfies
$N$-tuplewise independence.
\hfil\break  

   Remark 2.1(A) gives a standard scheme that was 
apparently first used by Doeblin [11] for the 
case $d = 1$ (random sequences),
in his construction of an i.i.d.\ sequence whose
marginal distribution is ``universal,'' i.e.\ in the 
domain of partial attraction to all infinitely divisible 
laws.
To verify statement (A), one can first observe with a
simple argument that
for any Borel sets $B_{n,k} \subset {\bf R}$,
$n \in {\bf N},\ k \in {\bf Z}^d$, and 
any $i \in {\bf Z}^d$,
$$P\Bigl(\bigcap_{(n,k) \in 
{\bf N} \times ({\bf Z} \uparrow d)}
\{Y^{(n)}_k \in B_{n,k}\}\Bigl) = 
P\Bigl(\bigcap_{(n,k) \in 
{\bf N} \times ({\bf Z} \uparrow d)}
\{Y^{(n)}_{k+i} \in B_{n,k}\}\Bigl).
$$
After that, statement (A) (that is, the strict
stationarity of $X$) follows from an elementary
measure-theoretic argument.
Then Property (B) follows by an elementary
argument.
\hfil\break

  {\bf Remark 2.2.} (A) Suppose $d$ is a positive
integer.  Suppose $Y := (Y_k, k \in {\bf Z}^d)$ is a
(not necessarily stationary) random field.
Suppose that for each $u \in {\bf Z}^d$,
$Y^{(u)} := (Y^{(u)}_j, j \in {\bf Z}^d)$ is a random
field whose distribution on (the Borel $\sigma$-field on)
${\bf R}^{{\bf Z} \uparrow d}$ is the same as that of $Y$.
Suppose these random fields $Y^{(u)},\ u \in {\bf Z}^d$
are independent of each other.
Suppose $f: {\bf R}^{{\bf Z} \uparrow d} \to {\bf R}$ 
is a Borel function.
Suppose $X := (X_k, k \in {\bf Z}^d)$ is the random 
field defined by 
$X_k := f((Y^{(k-j)}_j, j \in {\bf Z}^d))$ for
$k \in {\bf Z}^d$.
Then the random field $X$ is strictly stationary.

   (B) In the context of (A) above, if also $N \geq 2$
is an integer and the random field
$Y$ satisfies $N$-tuplewise
independence, then the random field $X$ satisfies
$N$-tuplewise independence.
\hfil\break

   The scheme in Remark 2.2(A) was employed by Olshen [20]
in the case $d=1$ (random sequences) to ``convert'' a
non-stationary random sequence with certain 
``ergodic-theoretic'' properties 
(trivial ``past'' and ``future'' tail $\sigma$-fields and
a nontrivial ``double'' tail $\sigma$-field) into a
strictly stationary sequence with those properties.
For $d \geq 2$, it was similarly employed in [2]
to construct a strictly stationary random field 
$X := (X_k, k \in {\bf Z}^d)$ that satisfies 
$\rho(2) = 0$ and whose ``tail $\sigma$-field''
is (modulo sets of probability 0) identical to 
$\sigma(X)$ itself.
To verify statement (A), one can first observe with a
simple argument (slightly cumbersome to write out
in more detail) that
for any Borel sets $B_{k,j} \subset {\bf R}$,
$k \in {\bf Z}^d,\ j \in {\bf Z}^d$, and any 
$i \in {\bf Z}^d$,
$$P\Bigl(\bigcap_{(k,j) \in 
({\bf Z} \uparrow d) \times ({\bf Z} \uparrow d)}
\{Y^{(k-j)}_j \in B_{k,j}\}\Bigl) = 
P\Bigl(\bigcap_{(k,j) \in 
({\bf Z} \uparrow d) \times ({\bf Z} \uparrow d)}
\{Y^{((k+i)-j)}_j \in B_{k,j}\}\Bigl).
$$
After that, statement (A) (the strict stationarity 
of $X$) follows from an elementary
measure-theoretic argument.
Then Property (B) follows by an elementary argument.
\hfil\break

   {\bf Remark 2.3.} (A) Suppose $d$ is a positive integer
and $X := (X_k, k \in {\bf Z}^d)$ is a strictly
stationary random field.
Let $F: {\bf R} \to [0,1]$ denote the (marginal) 
distribution function of the random variable 
$X_{\bf 0}$ ---
that is, $F(x) := P(X_{\bf 0} \leq x$). 
Let $G: {\bf R} \to [0,1]$ be the function defined by
$G(x) := \lim_{y \to x-} F(y)$. 
Let $F^{-1}: (0,1) \to {\bf R}$ denote the 
``generalized inverse'' function of $F$ --- that is,   
$F^{-1}(u) := \inf \{x \in {\bf R}: F(x) \geq u\}$. 
Let $V := (V_k, k \in {\bf Z}^d)$ be a random field of
independent, identically distributed random variables,
each uniformly distributed on the interval $[0,1]$,
with this random field $V$ being independent of the
random field $X$.  Define the random field
$U := (U_k, k \in {\bf Z}^d)$ as follows:
For each $k \in {\bf Z}^d$, 
$$ U_k := G(X_k) + V_k \cdot [F(X_k) - G(X_k)].
\eqno (2.1) $$
Then this random field $U$ has the following properties:  
   (1) The random field $U$ is strictly stationary.   
   (2)  The random variable $U_{\bf 0}$ is
uniformly distributed on the interval $[0,1]$. \break
   (3) For each $k \in {\bf Z}^d$, 
$X_k = F^{-1}(U_k)$ a.s.
   (4) For any two nonempty disjoint subsets
$S$ and $T$ of ${\bf Z}^d$, 
$$ \eqalignno{
\alpha(\sigma(U_k, k \in S), \sigma(U_k, k \in T)) &=
\alpha(\sigma(X_k, k \in S), \sigma(X_k, k \in T))\
\quad {\rm and} & (2.2) \cr
\rho(\sigma(U_k, k \in S), \sigma(U_k, k \in T)) &=
\rho(\sigma(X_k, k \in S), \sigma(X_k, k \in T)); & (2.3)      
}$$
and hence for each positive integer $n$,
(i) $\alpha(U,n) = \alpha(X,n)$,
(ii) $\rho(U,n) = \rho(X,n)$, \break
(iii) $\rho'(U,n) = \rho'(X,n)$, and  
(iv) $\rho^*(U,n) = \rho^*(X,n)$.

   (B) In the context of (A) above, if also $N \geq 2$
is an integer and the random field $X$ satisfies
$N$-tuplewise independence, then the random field
$U$ satisfies $N$-tuplewise independence.
\hfil\break

   The construction in (2.1) is classic.   
Mixing properties of this construction (\`a la the 
entire sentence (4) above)
essentially go back to Donald W.K.\ Andrews 
and Manfred Denker, independently of each other, in
unpublished manuscripts around the year 1982, both in the
context of random sequences (the case $d=1$).
(See e.g.\ the abstract in [10].)\ \ 
One reference for the proof in the case $d=1$ is
[4, v1, Theorem 6.8].  
The proof for general $d \geq 1$ is the same as for $d=1$.
(In particular, for (2.2) and (2.3), refer to (2.1)
and see
[4, v1, Theorem 6.2(I)(II) and Remark 6.3 (its third
paragraph)].)\ \ 
In Tone's proof of her result in [24] cited above,
part (A) was applied in the context of random fields
(general $d \geq 1$), with $X_{\bf 0}$ having an
absolutely continuous distribution.  
Part (B) is a trivial addendum, an immediate consequence
of (2.1) and the properties of the random field $V$.
\hfil\break

   In our application of the schemes above, the 
following lemma will play a key role. 
It is due to Cs\'aki and Fischer [8, Theorem 6.2]. 
For a generously detained presentation of its proof,
see e.g.\ [4, v1, Theorem 6.1].
\hfil\break

   {\bf Lemma 2.4.} {\sl Suppose 
${\cal A}_1, {\cal A}_2, {\cal A}_3, \dots$ and
${\cal B}_1, {\cal B}_2, {\cal B}_3, \dots$ are
$\sigma$-fields, and the $\sigma$-fields
${\cal A}_n \vee {\cal B}_n$, $n \in {\bf N}$ are
independent.  
Then
$$ \rho\Bigl(\bigvee_{n \in {\bf N}} {\cal A}_n, 
\bigvee_{n \in {\bf N}} {\cal B}_n \Bigl)\ =\
\sup_{n \in {\bf N}} \rho({\cal A}_n, {\cal B}_n). $$}

   The next definition (an old classic one), 
and the lemma right after it,
will play a key role in obtaining $N$-tuplewise
independence in the constructions for 
Theorems 1.4 and 1.5.
The redundancy in this definition may help
highlight its main features.
\hfil\break

   {\bf Definition 2.5.} For any integer $m \geq 3$,  
define the probability measures  $\nu^{(m)}_0$  and
$\nu^{(m)}_1$ on $\{-1,1\}^m$ (the set of all $m$-tuples
of $+1$'s and $-1$'s) as follows:  For each
$x := (x_1, x_2, \dots, x_m) \in \{-1,1\}^m$,
$$ \nu^{(m)}_0(x) := 1/2^m \quad {\rm and} \quad
\nu^{(m)}_1(x) := \cases{
1/2^{m-1} & if 
$x_1 \cdot x_2 \cdot \dots \cdot x_m = 1$ \cr
0 & if $x_1 \cdot x_2 \cdot \dots \cdot x_m = -1.$
\cr } $$
For any integer $m \geq 3$ and any 
$\theta \in (0,1)$, define the
probability measure $\nu^{(m)}_\theta$ on $\{-1,1\}^m$ by
$$ \nu^{(m)}_\theta := (1-\theta) \nu^{(m)}_0 + 
\theta \nu^{(m)}_1. \eqno (2.4) $$  
By simple arithmetic, for a given 
$x := (x_1, x_2, \dots, x_m) \in \{-1,1\}^m$, the
number $\nu^{(m)}_\theta (x)$ is equal to
$(1 + \theta)/2^m$ (resp.\ $(1 - \theta)/2^m$) if
$x_1 \cdot x_2 \cdot \dots \cdot x_m = 1$
(resp.\ if $x_1 \cdot x_2 \cdot \dots \cdot x_m = -1$). 
\hfil\break
 
   {\bf Lemma 2.6.} {\sl Suppose $m$ is an integer such
that $m \geq 3$, and $\theta \in [0,1]$.
Suppose $V := (V_1, V_2, \dots , V_m)$ is an 
$\{-1,1\}^m$-valued random vector whose distribution is
$\nu^{(m)}_\theta$.
Then the following statements hold:

   (A) For any permutation $\sigma$ of $\{1, 2, \dots, m\}$,
the distribution of the random vector 
$V^{(\sigma)} := (V_{\sigma(1)}, V_{\sigma(2)}, \dots,
V_{\sigma(m)})$ is $\nu^{(m)}_\theta$ as well.

   (B) For each $i \in \{1, 2, \dots, m\}$, one has that
$P(V_i = 1) = P(V_i = -1) = 1/2$.

   (C) Every $m-1$ of the random variables 
$V_1, V_2, \dots, V_m$ are independent.

   (D) If $S$ and $T$ are nonempty disjoint sets
   whose union is $\{1, 2, \dots, m\}$, then
$$ 
4 \alpha(\sigma(V_i, i \in S), \sigma(V_i, i \in T))
= \rho(\sigma(V_i, i \in S), \sigma(V_i, i \in T))
= \theta. $$}

   Property (A) in Lemma 2.6 is an immediate 
consequence of the sentence after (2.4).
Properties (B), (C), and (D) are elementary and 
were shown in [1, Lemma 3.2].
Because of property (A), Definition 2.5 can be
put in a more flexible form, as follows: 
\hfil\break

   {\bf Definition 2.7.}  For any finite set $S$
with $m := {\rm card}\thinspace S \geq 3$ and any 
$\theta \in [0,1]$, a collection 
$(V_i, i \in S)$ of $\{-1,1\}$-valued random variables 
is said to have the distribution $\nu^{(m)}_\theta$
if the distribution of the random vector 
$V := (V_{s(1)}, V_{s(2)}, \dots, V_{s(m)})$ 
is $\nu^{(m)}_\theta$ for any, hence every, ordered 
listing $s(1), s(2), \dots, s(m)$ of the elements 
(each exactly once) of $S$.
\hfil\break

\noindent {\bf 3. Proof of Theorem 1.4} \hfil\break

   The proof of Theorem 1.4 will start with the following
lemma:
\hfil\break
        
   {\bf Lemma 3.1.}  {\sl Suppose $d$, $N$, and $n$ are
positive integers such that $N \geq 2$.
Suppose $\theta \in [0,1]$.
Then there exists a strictly stationary,
$N$-tuplewise independent random field
$X := (X_k, k \in {\bf Z}^d)$
such that
(i) $4 \alpha(X,n) = \rho^*(X,1) = \theta$ and
(ii) $\rho^*(X,n+1) = 0$.}  
\hfil\break

Note that by (1.8), properties (i) and (ii) here imply 
that
\hfil\break   
(a) $\rho'(X,1) = \rho(X,1) = 4\alpha(X,1) = \rho^*(X,n) = \rho'(X,n) = \rho(X,n) = \theta$ and \hfil\break
(b) $\rho'(X,n+1) = \rho(X,n+1) = \alpha(X,n+1) = 0$.
\hfil\break
Such consequences of (1.8) will occur regularly
in the rest of this note.
\hfil\break 

   {\bf Proof.} The proof will be divided
into two ``steps.''
\hfil\break

   {\bf Step 1.} Let $M \geq 2$ be an
integer such that $M^d -1 \geq N$. 
Define the set
$\Lambda := \{0, n, 2n, \dots, (M-1)n\}^d$.
This set has cardinality $M^d$.
\hfil\break   
   
   Let $Y := (Y_k, k \in {\bf Z}^d)$ be a random field 
with the following properties:
(i) The random variables $Y_k, k \in \Lambda$ take
only the values $-1$ and $+1$; and the distribution of
this collection $(Y_k, k \in \Lambda)$ is
$\nu^{(M \uparrow d)}_\theta$ (see Definition 2.7).
(ii) The random variables $Y_k, k \in {\bf Z}^d - \Lambda$
are constant, defined by $Y_k := 0$.
\hfil\break

   By Lemma 2.6(C), every $M^d -1$ of the random
variables $Y_k, k \in {\Lambda}$ are independent.
It follows trivially that the entire random field
$Y$ satisfies $(M^d -1)$-tuplewise independence,
and hence (by the first sentence of Step 1 here)
$N$-tuplewise independence.
\bigskip  

   By Lemma 2.6(D), if $S$ and $T$ are nonempty,
disjoint subsets of $\Lambda$ and their union is 
$\Lambda$, then
$$4\alpha(\sigma(Y_k, k \in S), \sigma(Y_k, k \in T))
= \rho(\sigma(Y_k, k \in S), \sigma(Y_k, k \in T))
= \theta.$$
It follows from a simple argument that
$$ 4\alpha(Y,n) = \rho^*(Y,1) = \theta. 
\eqno (3.1) $$
(To see that $4\alpha(Y,n) = \theta$, refer to (1.4) and
consider the index sets $\{k \in {\bf Z}^d: k_1 \leq 0\}$ 
and $\{k \in {\bf Z}^d: k_1 \geq n\}$; each element of
$\Lambda$ belongs to one of those two sets, and each of
those two sets has at least one element of $\Lambda$.)
\hfil\break

   Also, if $S$ and $T$ are nonempty disjoint subsets
of ${\bf Z}^d$ such that 
${\rm dist}(S,T) \geq n+1$, and $S$ and $T$ each 
have at least one element of $\Lambda$, 
then some element of $\Lambda$
is in neither $S$ nor $T$.
(If $j \in \Lambda \cap S$ and $k \in \Lambda \cap T$,
then there exists a ``chain''
$(j := \kappa(0), \kappa(1), \kappa(2), \dots, 
\kappa(p) := k)$ 
of elements of $\Lambda$ such that
for each $i \in \{1, 2, \dots, p\}$, 
$\|\kappa(i) - \kappa(i-1)\| = n$ ---  with 
$\kappa(i)_u - \kappa(i-1)_u = \pm n$ for 
some $u \in \{1, \dots, d\}$;   
and if also $\Lambda \subset S \cup T$, then for some
$i \in \{1, 2, \dots, p\}$, one has
$\kappa(i-1) \in S$ and $\kappa(i) \in T$,
forcing ${\rm dist}(S,T) \leq n$, a contradiction.)\ \ 
It follows from Lemma 2.6(C) that
$$ \rho^*(Y,n+1) = 0.  \eqno (3.2) $$

   {\bf Step 2.} For each $u \in {\bf Z}^d$, let
$Y^{(u)} := (Y^{(u)}_j, j \in {\bf Z}^d)$ be a 
random field with the same distribution 
(on $\{-1, 0, 1\}^{{\bf Z} \uparrow d}$) 
as the random field $Y$ above.
Let these random fields $Y^{(u)}$, $u \in {\bf Z}^d$
be constructed in such a way that they are independent
of each other.
For technical convenience, assume that for {\it every\/}
$\omega \in \Omega$ and every $u \in {\bf Z}^d$,
$Y^{(u)}_j(\omega) \in \{-1,1\}$ (resp.\ $= 0$)
if $j \in \Lambda$ (resp.\ $j \in {\bf Z}^d - \Lambda$). 
\hfil\break

   Referring to the first paragraph of Step 1,
let $\phi$ be a one-to-one function from $\Lambda$ 
onto $\{ 0, 1, \dots, M^d -1\}$. 
Define the random field $X := (X_k, k \in {\bf Z}^d)$
as follows:  For each $k \in {\bf Z}^d$,
$$ X_k := \sum_{j \in \Lambda}
2^{\phi(j)} \cdot [(Y^{(k-j)}_j + 1)/2].  \eqno (3.3) $$
By Remark 2.2(A), 
this random field $X$ is strictly stationary;
and by Remark 2.2(B) and the third paragraph of Step 1 
above, it satisfies $N$-tuplewise independence.  
Our remaining task is to verify properties (i) and (ii)  
in the statement of Lemma 3.1. 
\hfil\break

   The quantity $(1+a)/2$ equals 0 (resp.\ 1) if
$a = -1$ (resp.\ 1).
Hence by the first paragraph of Step 2 here,
for a given $k \in {\bf Z}^d$, 
the random variables
$(Y^{(k-j)}_j + 1)/2$, $j \in \Lambda$
each take just the values 0 and 1.
Also, the binary expansion of a given positive
integer is unique.
It follows easily from (3.3) that for any given 
$k \in {\bf Z}^d$,
$$\sigma(X_k)  
= \sigma (Y^{(k-j)}_j, j \in {\Lambda})
= \sigma (Y^{(k-j)}_j, j \in {\bf Z}^d).  \eqno (3.4) $$
Here the second equality
follows from the fact that the random variables
$Y^{(k-j)}_j,\ j \in {\bf Z}^d - \Lambda$ are constant
(in fact 0).

   Now suppose $S$ and $T$ are nonempty  
subsets of ${\bf Z}^d$.  By (3.4),
$$ \sigma(X_k, k \in S) =
\bigvee_{j \in {\bf Z} \uparrow d}
\sigma(Y^{(j)}_\ell, \ell \in S-j), \eqno (3.5) $$
and the same holds with each $S$ replaced by $T$.
It follows trivially that
$$ \rho\Bigl(\sigma(X_k, k \in S), 
\sigma(X_k, k \in T)\Bigl)
\geq 
\rho\Bigl(\sigma(Y^{({\bf 0})}_k, k \in S), 
\sigma(Y^{({\bf 0})}_k, k \in T)\Bigl), \eqno (3.6)$$
and that the same holds with the symbol $\rho$
replaced on both sides by $\alpha$.  Also, by
(3.5) and its analog for $T$, and Lemma 2.4, one has that  
$$ \rho\Bigl(\sigma(X_k, k \in S), \sigma(X_k, k \in T)
\Bigl) =
\sup_{j \in {\bf Z} \uparrow d}
\rho\Bigl(\sigma(Y^{(j)}_\ell, \ell \in S-j),
\sigma(Y^{(j)}_\ell, \ell \in T-j)\Bigl). 
\indent \eqno (3.7) $$ 

If the sets $S$ and $T$ are disjoint, then for any
given $j \in {\bf Z}^d$, the sets $S-j$ and $T-j$
are disjoint.
Hence by (3.1) and (3.7),
$\rho^*(X,1) \leq \theta$.
However, by (3.6) and (3.1), one also has that
$\rho^*(X,1) \geq \rho^*(Y^{({\bf 0})},1) = \theta$.
Hence 
$$\rho^*(X,1) = \theta.  \eqno (3.8)$$
Now if the sets $S$ and $T$ satisfy 
${\rm dist}(S,T) \geq n+1$, then for each $j \in{\bf Z}^d$,
${\rm dist}(S-j,T-j) \geq n+1$, and hence 
by (3.7) and (3.2), 
$\rho(\sigma(X_k, k \in S), \sigma(X_k, k \in T)) = 0$.
It follows that 
$ \rho^*(X, n+1) = 0$.
Also, by (3.1) and the analog of (3.6) for 
$\alpha(\dots)$, one has that
$4\alpha(X,n) \geq 4\alpha(Y^{({\bf 0})}, n) = \theta$.
Hence by (1.8) and (3.8), $4\alpha(X,n) = \theta$.
All equations in properties (i) and (ii) in Lemma 3.1
have been verified, and the proof of Lemma 3.1 is
complete.
\hfil\break 

   {\bf Proof of Theorem 1.4.}  Suppose $d$ and $N$ and
the sequence $(c_1, c_2, c_3, \dots)$ are as in the
statement of Theorem 1.4.
For each positive integer $n$, applying Lemma 3.1, let 
$Z^{(n)} := (Z^{(n)}_k, k \in {\bf Z}^d)$ 
be a strictly stationary, $N$-tuplewise independent random field such that
$$4\alpha(Z^{(n)},n) = \rho^*(Z^{(n)},1) = c_n \quad 
{\rm and} \quad  \rho^*(Z^{(n)}, n+1) = 0. \eqno (3.9) $$
Let these random fields $Z^{(n)}, n \in {\bf N}$
be constructed in such a way that they are independent of
each other.

   Let $\psi: {\bf R} \times {\bf R} \times {\bf R} 
\times \dots \to {\bf R}$ 
be a function which is one-to-one, onto, and
bimeasurable (with respect to the Borel $\sigma$-fields).
(Such functions are well known to exist.)\ \
Define the random field $X := (X_k, k \in {\bf Z})$
as follows:  For each $k \in{\bf Z}^d$,
$$ X_k := \psi (Z^{(1)}_k, Z^{(2)}_k, Z^{(3)}_k, \dots). 
\eqno (3.10) $$

   By Remark 2.1(A)(B), the random field $X$ is strictly
stationary and satisfies $N$-tuplewise independence.
By (3.10) and the properties of the function $\psi$,
one has that for each element $k \in {\bf Z}^d$,
$\sigma(X_k) = 
\bigvee_{n \in {\bf N}}\sigma (Z^{(n)}_k)$.
Hence for each $n \in {\bf N}$, by (3.9),
$$ \alpha(X,n) \geq \alpha(Z^{(n)},n)
= c_n/4.  \eqno (3.11) $$
It also follows from Lemma 2.4 that for each $n \in {\bf N}$
$$ \rho^*(X,n) = 
\sup_{m \in {\bf N}} \rho^*(Z^{(m)}, n). \eqno (3.12) $$
Now for a given $m \in {\bf N}$ and $n \in {\bf N}$, 
one has that
$\rho^*(Z^{(m)}, n) = 0$ by (3.9) if ($n \geq 2$ and)
$m < n$; and if instead $m \geq n$ then
$\rho^*(Z^{(m)}, n) = c_m \leq c_n$ by (3.9), (1.8), and 
the hypothesis of Theorem 1.4.
Hence for each $n \in {\bf N}$, $\rho^*(X,n) = c_n$
by (3.12).
Hence for each $n \in {\bf N}$, by (3.11) and (1.8),
the equalities in property (C) of Theorem 1.4 all hold
for the random field $X$ here.

   All that remains is to convert the random field $X$
into one in which the (marginal) distribution of 
$X_{\bf 0}$ is uniform on $[0,1]$ without changing
the other properties stated in Theorem 1.4. 
On accomplishes that by applying Remark 2.3(A)(B).
That completes the proof of Theorem 1.4.
\hfil\break

\noindent {\bf 4. Proof of Theorem 1.5}
\hfil\break

   The proof of Theorem 1.5 will first involve two
lemmas.
The constructions in the proofs of those lemmas will
be somewhat related to the construction in [2].
\hfil\break

   {\bf Lemma 4.1.} {\sl Suppose $d$, $N$, and $n$ are positive integers such that $d \geq 2$ and $N \geq 2$.
Then there exists a strictly stationary,
$N$-tuplewise independent random field
$X := (X_k, k \in {\bf Z}^d)$ such that
(i) $\rho^*(X,n) =1$ and 
(ii) $\rho(X,1) = 1$ and $\rho(X,2) = 0$.}
\hfil\break

   {\bf Proof.} The proof will be divided
into two ``steps.''  
\hfil\break  
   
   {\bf Step 1.}  Increasing $n$ if necessary, 
we assume without loss of generality that 
$$  n > N.  \eqno (4.1) $$

   Let $\Gamma_0$ denote the set of all points
$k := (k_1, k_2, \dots, k_d) \in 
\{-n, -n+1, -n+2, \dots, n\}^d$
such that $k_u \in \{-n,n\}$ for at least one index 
$u \in \{1, 2, \dots, d\}$.  
That is, $\Gamma_0$ is the ``boundary'' or ``shell''
of the ``cube'' $\{-n, -n+1, \dots, n\}^d$. 
Define the set
$\Gamma := \Gamma_0 \cup \{{\bf 0}\}$.
Note that
$$ {\rm dist}\Bigl(\Gamma_0, \{{\bf 0}\}\Bigl) = n.  
\eqno (4.2) $$
Also, the sets $\{k \in {\bf Z}^d: k_1 \leq 0\}$ and
$\{k \in {\bf Z}^d: k_1 \geq 1\}$ each contain 
elements of $\Gamma$.
Also, for a given $j \in {\bf Z}$ and a given
index $u \in \{1, 2, \dots, d\}$, one of the following
three statements holds (depending on whether
$j \geq n+1$, $j \leq -n-1$, or $-n \leq j\leq n$): either
(a) $\Gamma \subset \{k \in {\bf Z}^d: k_u \leq j-1\}$, or
(b) $\Gamma \subset \{k \in {\bf Z}^d: k_u \geq j+1\}$, or
(c) the (``slice'') set
$\{k \in {\bf Z}^d: k_u = j\}$ contains elements of 
$\Gamma$ (here the assumption $d \geq 2$ is used). 
These trivial observations will be useful shortly.
\hfil\break

   Let $Y := (Y_k, k \in {\bf Z}^d)$ be a random field 
with the following properties:
(i) The random variables $Y_k, k \in \Gamma$ take
only the values $-1$ and $+1$; and the distribution of
this collection $(Y_k, k \in \Gamma)$ is
$\nu^{({\rm card}\thinspace \Gamma)}_1$
(see Definition 2.7 and note the ``extreme'' 
value 1 in the subscript here).
(ii) The random variables $Y_k, k \in {\bf Z}^d - \Gamma$
are constant, defined by $Y_k := 0$.
One has that
$$ \rho^*(Y,n) = 1 \quad {\rm and} \quad
\rho(Y,1) = 1 \quad {\rm and} \quad \rho(Y,2) = 0. 
\eqno (4.3) $$
In (4.3), the first equality holds by (4.2) and 
Lemma 2.6(D), and the second equality holds by 
Lemma 2.6(D) and the sentence right after (4.2). 
The third equality in (4.3) holds by the second sentence 
after (4.2). 
(In the definition of $\rho(Y,2)$, one can represent the
pairs of index sets \`a la (a) and (b) of the second
sentence after (4.2), and for $u,j$ such that  
(c) in that sentence holds, one applies Lemma 2.6(C)).
\hfil\break

By Lemma 2.6(C), every 
$({\rm card}\thinspace \Gamma) -1$ of the random
variables $Y_k, k \in {\Gamma}$ are independent.
Now trivially 
${\rm card}\thinspace \Gamma > n$; 
and now one has
by (4.1) that every $N$ of the random
variables $Y_k, k \in {\Gamma}$ are independent.  
It follows trivially that the entire random field
$Y$ satisfies $N$-tuplewise independence.
\hfil\break

   {\bf Step 2.}  Now we follow the argument in Step 2
of the proof of Lemma 3.1 in Section 3.
Here we shall just describe the changes.

   In the first paragraph of Step 2 there, replace each
$\Lambda$ by $\Gamma$.   
   
   In the second paragraph of Step 2 there, the changes 
are as follows:   
(a) One lets 
$\phi$ be a one-to-one function from $\Gamma$
onto $\{0,1,\dots, ({\rm card}\thinspace \Gamma) -1\}$.
(b) In eq.\ (3.3), the symbol $\Lambda$ is replaced by
$\Gamma$.
(c) The ``remaining task'' is to verify (i) and (ii) in
the statement of Lemma 4.1 (instead of in Lemma 3.1).

   In the third paragraph of Step 2 (including in
eq.\ (3.4)), the symbol $\Lambda$ is replaced throughout
by $\Gamma$.

   The fourth paragraph of Step 2 (with eqs.\ (3.5), (3.6), and (3.7)) remains unchanged.
The final (i.e.\ fifth) paragraph of Step 2 is modified
as follows:
First, one simply uses (4.3) and (3.6) (and (1.8)) to
obtain that $\rho^*(X,n) = 1$ and $\rho(X,1) = 1$.
Then one observes that if $S$ and $T$ are nonempty
subsets of ${\bf Z}^d$ that are on opposite sides of 
some ``slice'' (a set of the form 
$\{k \in {\bf Z}^d: k_u = h\}$ where $h \in {\bf Z}^d$
and $u \in \{1, \dots, d\}$), then for any 
$j \in {\bf Z}^d$, that is true as
well for the sets $S-j$ and $T-j$ (with respect to
a possibly different ``slice'').
From that and (3.7) and the third equality in (4.3),
one obtains $\rho(X,2) = 0$.
Thus properties (i) and (ii) in Lemma 4.1 have been
verified, and the proof of that lemma is complete.
\hfil\break

   {\bf Lemma 4.2.} {\sl Suppose $d$ and $N$ are 
integers such that $d \geq 2$ and $N \geq 2$.
Then there exists a strictly stationary,
$N$-tuplewise independent random field
$X := (X_k, k \in {\bf Z}^d)$ such that
(i) $\rho^*(X,n) =1$ for all $n \in {\bf N}$ and 
(ii) $\rho(X,1) = 1$ and $\rho(X,2) = \rho'(X,2)= 0$.}
\hfil\break

   {\bf Proof.}  One follows the Proof of Theorem 1.4
(after the proof Lemma 3.1) in Section 3,
but with the following changes:

   In the first paragraph of that proof, the mention
of a sequence $(c_1, c_2, c_3, \dots)$ is to be omitted,
one applies Lemma 4.1 (instead of Lemma 3.1), and
eq.\ (3.9) is replaced by
$$ \rho^*(Z^{(n)},n) = 1 \quad {\rm and} \quad
\rho(Z^{(n)},1) = 1 \quad {\rm and} \quad
\rho(Z^{(n)},2) = 0. \eqno (4.4) $$

   Eq.\ (3.10) and its entire paragraph, and also the 
first two sentences after (3.10), are left unchanged.
From (3.10) (which yields 
$\rho^*(X,n) \geq \rho^*(Z^{(n)},n)$
for each $n \in {\bf N}$) and (4.4), one
obtains that $\rho^*(X,n) = 1$ for all $n \in {\bf N}$.
Similarly from (3.10) and (4.4) one trivially obtains
$\rho(X,1) = 1$; and by (3.10), (4.4), and (say) Lemma 2.4,
one also obtains $\rho(X,2) = 0$.
From that last equality, one also has that 
$\rho'(X,2) = 0$, by an elementary argument (see e.g.\ 
[4, v3, Proposition 29.5]).  
That completes the proof of Lemma 4.2.
\hfil\break

   {\bf Proof of Theorem 1.5.} Suppose the integers
$d$ and $N$ and the sequence $(c_1, c_2, c_3, 
\allowbreak \dots)$
are as in the statement of Theorem 1.5.
Let $Y := (Y_k, k \in {\bf Z}^d)$ be a random field that
satisfies all properties stated for the random field
$X$ in Theorem 1.4.
Let $Z := (Z_k, k \in {\bf Z}^d)$ be a random field that
satisfies all properties stated for the random field
$X$ in Lemma 4.2.
Let these two random fields $Y$ and $Z$ be constructed
in such as way that they are independent of each other.
Let $\zeta : {\bf R} \times {\bf R} \to {\bf R}$ 
be a function which is one-to-one, onto, and 
bimeasurable (with respect to the Borel $\sigma$-fields).
Define the random field $X := (X_k, k \in {\bf Z}^d)$ as
follows:  For each $k \in {\bf Z}^d$,
$$  X_k := \zeta(Y_k, Z_k). \eqno (4.5) $$

By Remark 2.1(A)(B) (applied to $Y$ and $Z$ and, say, a
sequence of degenerate random fields whose random
variables take only the value 0), the random
field $X$ is strictly stationary and $N$-tuplewise
independent.
By (4.5) and the properties of the function $\zeta$,
one has that for any given $k \in {\bf Z}^d$,
$$ \sigma(X_k) = \sigma(Y_k, Z_k). \eqno (4.6) $$
From (4.6) and Lemma 4.2 (for the random field $Z$), one
immediately obtains that
$\rho^*(X,n) \allowbreak = 1$ for all $n \in {\bf N}$, 
and that
$\rho(X,1) = 1$, and hence also $\rho'(X,1) = 1$.
Also, for each $n \geq 2$, one has by Lemma 2.4, 
followed by the properties in Theorem 1.4 (for $Y$) 
and Lemma 4.2 (for $Z$),
$$ \rho(X,n) = \max \{ \rho(Y,n), \rho(Z,n) \} = 
\max\{c_n,0\} = c_n, $$
and the same holds with each symbol $\rho$ replaced by
$\rho'$.

   Now all that remains is 
to convert the random field $X$ into one in which the (marginal) distribution of $X_{\bf 0}$ is uniform on 
$[0,1]$ without changing the other properties stated
in Theorem 1.5.
One accomplishes that with an application of 
Remark 2.3(A)(B).
That completes the proof of Theorem 1.5.
\hfil\break

\centerline {\bf References}
\bigskip

\refs [1] R.C.\ Bradley. 
Some examples of mixing random fields.
{\it Rocky Mountain J.\ Math.\/}
23 (1993) 495-519.

\refs [2] R.C.\ Bradley.  On regularity conditions for
random fields.  {\it Proc.\ Amer.\ Math.\ Soc.\/} 
121 (1994) 593-598.

\refs [3] R.C.\ Bradley.  On mixing rates for
nonisotropic random fields.
{\it Acta Sci.\ Math.\ (Szeged)\/} 65 (1999) 749-765. 

\refs [4] R.C.\ Bradley. 
{\it Introduction to Strong Mixing Conditions\/},
Volumes 1, 2, and 3. \allowbreak
Kendrick Press, Heber City (Utah), 2007.

\refs [5] R.C.\ Bradley. 
On the dependence coefficients associated with three 
mixing conditions for random fields,
In:\ {\it Dependence in Probability, Analysis and Number
Theory\/}, 
(I.\ Berkes, R.C.\ Bradley, H.\ Dehling,
M.\ Peligrad, and R.\ Tichy, eds.), 
pp.\ 89-121.
Kendrick Press, Heber City (Utah), 2010.

\refs [6] R.C.\ Bradley and A.R.\ Pruss.
A strictly stationary, $N$-tuplewise independent
counterexample to the central limit theorem.
{\it Stochastic Process.\ Appl.\/} 119 (2009) 
3300-3318.

\refs [7] R.C.\ Bradley and L.T.\ Tran.
Density estimation for nonisotropic random fields.
{\it J.\ Statist.\ Plann.\ Inference\/} 81 (1999) 51-70. 

\refs [8] P.\ Cs\'aki and J.\ Fischer. 
On the general notion of maximal correlation.
{\it Magyar Tud.\ Akad.\ Mat.\ Kutato Int.\ Kozl.\/}
8 (1963) 27-51.

\refs [9] A. Cuesta and C.\ Matr\'an.
On the asymptotic behavior of sums of pairwise
independent random variables.
{\it Statist.\ Probab.\ Lett.\/} 11 (1991) 201-210.

\refs [10]  M.\ Denker.  
Substituting independent processes.
In: {\it Journees de Theorie Ergodique C.I.R.M.\/}, 
Luminy, Marseille, July 5-10, 1982.

\refs [11] W.\ Doeblin.  
Sur l'ensemble de puissances d'une loi de probabilit\'e.  
{\it Studia Math.\/} 9 (1940) 71-96.

\refs [12] N.\ Etemadi.
An elementary proof of the strong law of large numbers.
{\it Z.\ Wahrsch.\ verw.\ Gebiete\/} 55 (1981) 119-122.

\refs [13] N.\ Etemadi. 
Convergence of weighted averages of random variables revisited.  
{\it Proc.\ Amer.\ Math.\ Soc.\/} 134 (2006) 2739-2744.

\refs [14] N.\ Etemadi and A.\ Lenzhen. 
Convergence of sequences of pairwise independent 
random variables.  
{\it Proc.\ Amer.\ Math.\ Soc.\/} 132 (2004) 1201-1202.

\refs [15] C.M\ Goldie and P.E\ Greenwood.  Variance of
set-indexed sums of mixing random variables and weak convergence of set-indexed processes.
{\it Ann.\ Probab.\/} 14 (1986) 817-839.

\refs [16] V.V.\ Gorodetski\u i.  The central limit theorem 
and an invariance principle for weakly dependent random fields. {\it Soviet Math.\ Dokl.\/} 29 (1984) 529-532.

\refs [17] S.\ Janson.  Some pairwise independent
sequences for which the central limit theorem fails.
{\it Stochastics\/} 23 (1988) 439-448.

\refs [18] H.\ Kesten and G.L.\ O'Brien.  Examples of
mixing sequences.  {\it Duke Math.\ J.\/} 43 (1976) 
405-415.

\refs [19] F.S.\ McWilliams and N.J.A.\ Sloane.
{\it The Theory of Error-Correcting Codes\/}.
North-Holland, Amsterdam, 1977.

\refs [20] R.A.\ Olshen.  The coincidence of measure
algebras under an exchangeable probability.  
{\it Z.\ Wahrsch.\ verw.\ Gebiete\/}
18 (1971) 153-158.

\refs [21] M.\ Peligrad and A.\ Gut [1999].  Almost sure
results for a class of dependent random variables.
{\it J.\ Theoret.\ Probab.\/} 
12 (1999) 87-104. 

\refs [22] M.\ Rosenblatt.  {\it Stationary Sequences 
and Random Fields.\/} Birkh\"auser, Boston, 1985.

\refs [23] M.\ Rosenblatt and D.\ Slepian. 
$n$th order Markov chains with every $N$ variables independent.  
{\it J.\ Soc.\ Indust.\ Appl.\ Math.\/} 
10 (1962) 537-549.

\refs [24] C.\ Tone.  Central limit theorems for strictly stationary random fields under strong mixing conditions.
Ph.D.\ Dissertation, Indiana University, 
Bloomington, 2010.

\refs [25] C.\ Tone.
Central limit theorems for Hilbert-space valued random 
fields satisfying a strong mixing condition. 
{\it Lat.\ Am.\ J.\ Probab.\ Math.\ Stat.\/} 
8 (2011) 77-94.

\bye